\documentclass{amsart}
\usepackage{amsfonts,amsmath,amsthm}
\usepackage{amssymb,latexsym} 
\newtheorem{theorem}{Theorem}[section] 
\newtheorem{lemma}[theorem]{Lemma}
\theoremstyle{definition} 
\newtheorem{definition}{Definition}[section] 
\theoremstyle{remark} 
\numberwithin{equation}{section} 
\newcommand{\real}{\mathbb{R}}

\newcommand{\complex}{\mathbb{C}}

\newcommand{\loc}{{\scriptstyle{loc}}} 
\newcommand{\ca}{\mathcal} 
\DeclareMathOperator{\supp}{Supp}

\newcommand{\ddt}{\frac{d \hfill}{dt}} 
\newcommand{\grad}{\nabla}
\newcommand{\eqdef}{\mathrel{:=}}

\newcommand{\RR}{\real}
\newcommand{\lpp}{\left(}
\newcommand{\rpp}{\right)}
\newcommand{\lb}{\left[}
\newcommand{\rb}{\right]}
\newcommand{\lB}{\left\{}
\newcommand{\rB}{\right\}}
\newcommand{\fakeq}{\phantom{=\;}}
\newcommand{\eval}[1]{{\big|_{#1}}}
\newcommand{\Eval}[2]{\mathchoice{\bigg|_{#1}^{#2}}{\big|_{#1}^{#2}}%
{\big|_{#1}^{#2}}{\big|_{#1}^{#2}}}
\begin{document} 

\title[Equivalence of Euler and Birkhoff-Rott]{A Criterion for the Equivalence of the 
Birkhoff-Rott and Euler Descriptions of Vortex Sheet Evolution} 

\author[M. Lopes Filho, H. Nussenzveig Lopes and S. Schochet]{Milton C. Lopes Filho}  
\address{Depto. de Matem\'{a}tica, IMECC-UNICAMP,
Cx. Postal 6065, Campinas SP 13081-970, Brazil.
e-mail: mlopes@ime.unicamp.br} 

\author[]{Helena J. Nussenzveig Lopes}  
\address{Depto. de Matem\'{a}tica, IMECC-UNICAMP,
Cx. Postal 6065, Campinas SP 13081-970, Brazil.
e-mail: hlopes@ime.unicamp.br}
 
\author[]{Steven Schochet}  
\address{School of Mathematical Sciences, Tel Aviv University, 
Ramat Aviv, Tel Aviv 69978 Israel.
e-mail: schochet@post.tau.ac.il} 

\begin{abstract}
  In this article we consider the evolution of vortex sheets in the plane
  both as a weak solution of the two dimensional incompressible Euler
  equations and as a (weak) solution of the Birkhoff-Rott equations.
  We begin by discussing the classical Birkhoff-Rott equations with 
  respect to arbitrary parametrizations of the sheet. We introduce a 
  notion of weak solution to the Birkhoff-Rott system and we prove 
  consistency of this notion with the classical formulation of the equations. 
  Our main purpose in this paper is to present a sharp criterion for the equivalence of 
  the weak Euler and weak Birkhoff-Rott descriptions of vortex sheet dynamics.
\end{abstract}

\maketitle


\section{Introduction}

There are two distinct points of view in the mathematical description 
of interface dynamics. The more natural one is to explicitly propagate the
interface itself using a time-dependent parametrization. 
An alternative approach is to embed the interface into the solution 
of a partial differential equation which can be evolved, carrying the
interface with it. We will refer to the former as the {\it explicit} approach
to interface dynamics, while the latter will be called the {\it implicit} 
approach. See \cite{of03} for a broad discussion and several instances of this 
dychotomy.

Vortex sheet evolution in two-dimensional, incompressible, ideal fluid
flow is a classical example of interface dynamics for which both
points of view have been widely addressed. The explicit approach in this 
context makes use of the Birkhoff-Rott equations. This system was originally 
derived by G. D. Birkhoff \cite{birkhoff} and it is implicit in the work of N. Rott 
\cite{rott56}. For the implicit point of view one uses the incompressible 
2D Euler equations, regarding the vortex sheet as a 
feature of a suitably defined weak solution.  The purpose of the present work 
is to establish a sharp smoothness criterion for the equivalence of these 
descriptions of vortex sheet motion. 

The physically meaningful notion of generalized solution for the Euler 
equations is the weak form of the momentum equations for velocity and pressure, known
as the weak velocity formulation. 
Indeed, the weak velocity formulation is an integral form of conservation of momentum, and 
therefore, it is close to physical first principles. In the problem of vortex sheet 
evolution, an alternative weak formulation, known as the weak vorticity 
formulation, has proven to be more tractable.  The weak vorticity formulation 
has been shown to be equivalent 
to the weak velocity formulation in all situations under consideration in this work, 
see \cite{delort,S}. In this context, it is natural to ask  
whether the Birkhoff-Rott description of vortex sheet evolution is  equivalent 
to these weak formulations as well. 

The problem of equivalence between Birkhoff-Rott and the weak forms of 2D Euler
may be regarded as solved in the case of vortex sheets with smooth densities on
smooth curves. Although there is no complete proof available, the proof of 
Theorem 6.1.2 in \cite{MP}  (which assumes that the vortex sheet is a graph) 
can be easily adapted to establish such a result in general. On the other hand, 
it was recently shown, see \cite{lopes1}, that a well-known exact solution of the 
Birkhoff-Rott equations fails to satisfy both the weak vorticity and the weak 
velocity formulations. Our main result is to show that, if the vortex sheet is a 
{\emph{regular}} curve, for each fixed time, then the condition for equivalence 
is that the vorticity density be square-integrable with respect to arclength.
The example in \cite{lopes1} shows that the integrability condition on the
vorticity density is sharp. The definition of regular curve is standard in
harmonic analysis, see \cite{david}; a rectifiable curve is \emph{regular} 
if its intersection with a ball of radius $r$ has length $\mathcal{O}(r)$. 
 
In this work, we require a weak formulation of the Birkhoff-Rott equations. 
This is not a standard topic, and therefore, we must treat it at length.
A special case of the weak Birkhoff-Rott equations formulated here was
given in equation (6.1.14) of \cite{MP} for vortex sheets which are graphs. 
Furthermore, the argument used in \cite{MP} to establish equivalence between
Birkhoff-Rott and 2D Euler in the smooth case should, in principle, be 
extendable to the weak solutions considered here. The method of proof we use is 
completely different. To motivate our particular weak formulation of the
Birkhoff-Rott equations we have stated and proved a consistency result, namely, 
that the weak formulation of Birkhoff-Rott plus minimal regularity assumptions is 
equivalent to classical Birkhoff-Rott.   

In order to properly contextualize our results we will give a brief account of 
the literature on vortex sheet evolution. Vortex sheets are a classical
topic in fluid dynamics. The complicated
evolution of vortex sheets is a natural source for the spontaneous
appearance of small scale motion in incompressible fluids, an
observation dating back to H. Helmholtz in 1868 (see the discussion
and references in \cite{batchelor}). This motivates the continuing
interest of the topic. The source of the small scales can be 
identified with a feature of vortex sheet motion, 
known as {\it Kelvin-Helmholtz} instability, see \cite{saffman} and 
references therein. 

In the classical paper \cite{moore} D. Moore presented a
theoretical account of how the Kelvin-Helmholtz instability could
drive a real analytic vortex sheet to form curvature singularities 
in finite time.  This work was later made rigorous by a number of authors, 
including results on existence of a 
local-in-time solution to Birkhoff-Rott for analytic initial data,
\cite{caflischo,duchonr,bardosfss} and formation of Moore's singularity
\cite{caflischo2,duchonr}. After singularity formation one
expects vortex sheet roll-up. Numerical studies, see for example \cite{krasny2},
illustrate the expected presence of double-branched spiral vortex sheets after
singularity formation for periodic perturbations of a planar sheet.
One is naturally led to study self-similar spiral shaped vortex sheets. 
The existence of exact self-similar spiral solutions of the Birkhoff-Rott 
equations goes back to an example due to L.Prandtl, generalized by T. Kambe,
see \cite{kambe,prandtl22}. These are examples of finite length 
logarithmic spirals. However, the generic roll-up after Moore's singularity 
and the roll-up of the elliptically loaded wing (another classical example),
\cite{krasny86desing,krasny2,krasny87,bakero} seem to lead to {\it infinite}
length algebraic spirals for which no rigorous existence is known. An
asymptotic description of self-similar solutions of the Birkhoff-Rott
equations with algebraic spirals was first presented by
Kaden~\cite{kaden31} and generalized by Pullin~\cite{pullins}.
The Birkhoff-Rott equations are elliptic in nature, and there
is a strong analogy between the Kelvin-Helmholtz instability and the 
Hadamard instability of the Cauchy problem for Laplace's equation.
In particular, the explicit approach to the vortex sheet evolution problem
is rather ill posed.  Rigorous results in this direction have recently 
been presented by G. Lebeau and also by S. Wu, see \cite{lebeau,wu}. 

The implicit approach to vortex sheets was pioneered by R. DiPerna and A. Majda
in a series of papers \cite{dipernam1,dipernam2,dipernam3}, where they
outlined a program for proving existence of weak solutions for the incompressible
2D Euler equations with vortex sheets as initial data. The DiPerna-Majda program 
was carried out in the case of steady weak solutions but remains open in
general~\cite{evans,majda1,NL,zheng}. In 1990, J.-M.
Delort proved the existence of a global-in-time weak solution for the
vortex sheet initial data problem with {\it distinguished sign}~\cite{delort},  
see also \cite{evansm,liux,majda2,S,S2}. A global-in-time existence result has 
also been proved in the case of mirror-symmetric flows with distinguished 
sign vorticity on each side of the mirror~\cite{lopes3}. Delort's Theorem and 
its mirror-symmetric extension provide the existence of a meaningful evolution 
for certain vortex sheets beyond singularity formation but give no information 
on their structure.  

The remainder of this article is divided into three sections. In the next
section we describe various formulations of both the Euler equations and the
Birkhoff-Rott equations, and we discuss the consistency of the weak
formulation of the latter.  The following section contains the precise
statement and proof of the equivalence between weak Birkhoff-Rott and
the weak vorticity formulation of 2D Euler, as well as the discussion
of the sharpness of this result. The final section contains some
interpretations of the work presented, as well as open problems and
conclusions. This work contains a complete answer to a problem formulated
by S. Wu at the end of \cite{wu}.

\section{Vortex sheet equations in weak form}
\label{sec:form}

This section is divided in three subsections. In the first we recall 
the weak vorticity formulation of the 2D incompressible Euler 
equations, implementing it in the special case where vorticity is 
concentrated on a curve. In the second subsection we discuss the
derivation of the several forms of the Birkhoff-Rott system.
In the third we introduce a notion of
weak solution of the Birkhoff-Rott equations and we prove consistency
of this notion with the classical form of the equations.

\subsection{Weak forms for the vorticity equation}
The vorticity formulation of the 2D incompressible Euler equations is
\begin{align} 
\label{eq:vort}
\omega_t + u \cdot \nabla \omega &= 0,
\\
\label{eq:bsvart} 
u &= K \ast \omega, 
\\ 
\label{eq:initome}
\omega(x,0) &= \omega_0(x),
\end{align}
with
\begin{equation} \label{K} K(x) = \frac{x^{\perp}}{2\pi |x|^2}, \end{equation}
where $(x_1,x_2)^{\perp} = (-x_2,x_1)$, $\omega$ is the vorticity, $u = (u_1,u_2)$ is the velocity and $K =
(K_1,K_2)$ is the kernel of the Biot-Savart law.
This system of equations can be reformulated
in several different ways that are relevant for the comparison with
the Birkhoff-Rott equations. First, the Lagrangian representation is
obtained by noting that \eqref{eq:vort} says that $\omega$ is preserved
along the particle trajectories having velocity $u$, i.e.,
\begin{align}
  \label{eq:lagrep}
    \omega(\Psi(x_0,t),t)&=\omega_0(x_0), \qquad\text{ where }
\\
\label{eq:vortraj}
\ddt\Psi(x_0,t)&=u(\Psi(x_0,t),t), 
\\
\label{eq:phinit}
\Psi(x_0,0)&=x_0, 
  \end{align}
  with $u$ defined as before by the Biot-Savart law~\eqref{eq:bsvart}.
  
  Vortex sheet flows have vorticities which are Radon measures supported on 
  rectifiable curves. In order to study such flows we require a weak formulation of the vorticity equation. 
  There are actually two versions: The traditional weak vorticity formulation
  (e.g. \cite{MP}) is obtained by multiplying \eqref{eq:vort} by
  a smooth test function~$\varphi$ having compact support in
  $[0,T)\times\RR^2$, integrating over space and time, and integrating
  by parts, which yields
\begin{equation}
  \label{eq:twvort}
  \begin{aligned}
    \int_0^T \int_{\RR^2}\omega(x,t)
\lB\varphi_t +u(x,t)\cdot\grad\varphi\rB\,dx\,dt
\ +\ \int_{\RR^2}\omega_0(x)\varphi(x,0)\,dx=0,
  \end{aligned}
\end{equation}
with $u$ still given by \eqref{eq:bsvart}. 

The velocity associated with a vortex sheet is, in principle, discontinuous
on the sheet, a fact which we will discuss in great detail later. 
The discontinuity at the sheet implies that the term 
$\omega u \cdot \nabla\varphi$ appearing in \eqref{eq:twvort}
is not well defined. 
However, this difficulty can be overcome by considering an
alternative weak formulation as follows.

The modern weak form can be
obtained \cite{S} from \eqref{eq:twvort} by substituting $u$ by
$K\ast\omega$, see \eqref{eq:bsvart}, and replacing the factor
multiplying $\omega(x,t)\omega(y,t)$ in the result by its symmetric
part. This yields
\begin{multline} \label{wvfid}
\int_0^{T} \int_{\real^2} \varphi_t \omega(x,t)\,dx\,dt +   
  \int_0^T \int_{\real^2} 
\int_{\real^2} H_{\varphi}(x,y,t)\omega(y,t)\omega(x,t) \,dx\,dy\,dt 
\\  + \int_{\real^2} \varphi(x,0)\omega_0(x)\,dx = 0, 
\end{multline}
where the function
\begin{equation*}
H_{\varphi}(x,y,t) \eqdef
\frac{\nabla\varphi(x,t)-\nabla\varphi(y,t)}{2} \cdot K(x-y)
\end{equation*}
is continuous for $x\ne y$, and it is also
bounded.  

Let $\mathcal{CBM}$ denote the space of bounded Radon measures
with no discrete part. 

\begin{definition} \label{wvf}
  Let $\omega \in L^{\infty}((0,T);\mathcal{CBM}(\real^2))$.
  We will say that $\omega$ is a weak solution of the
  Euler equations with initial data $\omega_0$ if, for any test 
  function $\varphi \in C^{\infty}_c([0,T)
  \times \real^2)$, \eqref{wvfid} holds.
\end{definition}

For smooth vorticities $\omega$ decaying sufficient rapidly at
infinity, all four formulations \eqref{eq:vort}--\eqref{eq:initome},
\eqref{eq:lagrep}--\eqref{eq:phinit} plus \eqref{eq:bsvart},
\eqref{eq:twvort} plus \eqref{eq:bsvart}, and \eqref{wvfid} are
equivalent \cite{delort,MP,saffman,S} to each
other and to both the classical and weak \cite{dipernam1} velocity
formulations of the Euler equations.  Moreover
\cite{delort,S}, the modern weak formulation remains
equivalent to the weak velocity formulation assuming only
that the vorticity lies in $L^\infty([0,T],\ca{BM})\cap
L^\infty([0,T],H^{-1}_{\loc})$.

It was shown in \cite{S} that the modern weak vorticity formulation makes sense
whenever the vorticity is a bounded measure having no discrete part. 
In particular, it makes sense when the vorticity is a measure concentrated 
along a smooth time-dependent curve. 

Let $\ca{C}_t$ be a smooth, time-dependent curve,
\begin{equation}\label{eq:ct}
\ca{C}_t\eqdef\{\xi = \xi(s,t)\mid s_0\le s\le s_1\}
\end{equation}
parametrized by arclength. Let $\gamma= \gamma(s,t)$ be a smooth density
and specify the vorticity to be the measure $\omega = \omega(x,t) = 
\gamma \delta_{\ca{C}_t}$ defined through the identity
\begin{equation}
  \label{eq:sheetvort}
  \langle \omega(\cdot,t), \varphi \rangle  \equiv \int_{s_0}^{s_1}
\gamma(s,t) \varphi(\xi(s,t))ds,
\end{equation}
for any test function $\varphi \in C^{\infty}_c(\real^2)$.
We assume that the initial vorticity is of the same form. 
If $\omega$ satisfies Definition \ref{wvf} then we substitute 
\eqref{eq:sheetvort} into \eqref{wvfid} to get
\begin{equation} \label{vsws}
\begin{array} {l}
 \displaystyle{\int_0^{T} \int_{s_0(t)}^{s_1(t)}} 
\varphi_t(\xi(s,t),t) \gamma(s,t) \, dsdt  \\ \\
 + \displaystyle{ \int_0^{T} \int_{s_0(t)}^{s_1(t)} 
\int_{s_0(t)}^{s_1(t)}} 
H_{\varphi}(\xi(r,t),\xi(s,t),t)\gamma(r,t)\gamma(s,t) \, drdsdt  \\ \\
+ \displaystyle{\int_{s_0(0)}^{s_1(0)}} \varphi(\xi_0(s),0)\gamma_0(s)\, ds = 0, 
\end{array} 
\end{equation}
for any test function $\varphi \in C^{\infty}_c([0,T) \times \real^2)$.

 
\subsection{The Birkhoff-Rott system}

We will now change our point of view, discussing  the  explicit approach
to vortex sheet evolution. Our objective is to examine the equation for 
the evolution of the sheet with respect to an arbitrary parametrization. 
To this end let us begin by considering the {\it linear} problem of 
transport, by a smooth vector field, of a measure concentrated on a smooth 
curve in the plane. 

Fix $v$ a given smooth vector field in the plane. Denote by $X=X_t$
the flow to time $t$ generated by $v$. Let $\mu_0$ be a Radon measure 
on the plane. We say that $\mu = \mu(\cdot, t)$ is the transport by 
$v$ of the measure $\mu_0$ if, for any Borelian subset 
$E \subseteq \real^2$, we have:
\[\mu(E,t) = \mu_0(X_{-t}(E)).\]
It is not hard to see that, if $\mu_0 = \gamma_0 \delta_{\ca{C}_0}$, 
then the transport by $v$ is of the form $\mu(\cdot,t) = 
\gamma \delta_{\ca{C}_t}$ with $\ca{C}_t = X_t(\ca{C}_0)$. 
Furthermore, under the same hypothesis', 
if $\Sigma_t \subseteq \ca{C}_t$ is the transport of a portion 
$\Sigma_0$ of $\ca{C}_0$, then 
\begin{equation} \label{dell}
 \frac{d}{dt} \int_{\Sigma_t} \gamma \cdot \, d\ell = 0.
\end{equation}

\begin{lemma} \label{lemmachtung}
Let $\mu_0 = \gamma_0\delta_{\ca{C}_0}$ be a Radon measure with support on a 
smooth curve $\ca{C}_0$. Let $\mu = \mu(\cdot,t)$ be the transport by 
$v$ of the measure $\mu_0$.

Let $y = y(\eta,t)$ be a parametrization of $\ca{C}_t$ and denote 
\[\sigma(\eta,t) \eqdef  \gamma(s(\eta,t),t) \frac{\partial s}{\partial \eta},\]
where $s=s(\eta,t)$ is arclength with respect to a reference point. 

Then there exists $a=a(\eta,t)$ such that the following equations are satisfied:
\begin{equation}  \label{achtung}
\left\{ \begin{array}{l}
 y_t + a y_{\eta} = v\\
 \sigma_t + (a\sigma)_{\eta}=0.\end{array} \right.\end{equation}
\end{lemma}
 
\begin{proof}
Consider a parametrization of $\ca{C}_0$, $z_0=z_0(\theta)$, $\theta \in I \subseteq \real$ 
and we assume for convenience that $0 \in I$. Let 
$z(\theta,t)=X_t(z_0(\theta))$ be a (Lagrangian) parametrization of $\ca{C}_t$. Let 
$s = s(\theta,t)$ be the arclength 
along $\ca{C}_t$ between $z(0,t)$ and $z(\theta,t)$. Then $\theta \mapsto s(\theta,t)$
is an invertible change-of-variables. We write $\xi = \xi(s,t)$ for the parametrization
with respect to the arclength $s$, measured from $z(0,t)$. 
 
With this notation it is a straightforward calculation to verify that 
\[\xi_t + \widetilde{a}(s,t)\xi_s  = v,\]
with 
\[\widetilde{a}(s,t) = \frac{\partial s}{\partial t} (\theta(s,t),t). \]

Next, implementing the condition \eqref{dell} gives, for any $\theta_0$, $\theta_1$,
\[ \frac{d}{dt}\int_{s(\theta_0,t)}^{s(\theta_1,t)} \gamma(s,t) ds = 0.\]
From this integral equation it follows easily that
\[\gamma_t + (\widetilde{a}(s,t)\gamma)_s = 0.\]

If $y=y(\eta,t)$ is any other parametrization of $\ca{C}_t$ and if 
$\sigma = \gamma \, s_{\eta}$ then it is immediate that
$y$ and $\sigma$ satisfy \eqref{achtung} with 
\[a(\eta,t) = \frac{\widetilde{a}(s(\eta,t),t) - s_t(\eta,t)}{s_{\eta}(\eta,t)}.\]
This concludes the proof.

\end{proof}

{\bf Remark:} System \eqref{achtung} is an explicit description of the
propagation of a curve which corresponds to the implicit description given
by the equation $\mu_t + \mbox{ div }(v\mu) = 0$, in the sense of 
distributions. Note that the function $a$ is a free parameter in \eqref{achtung},
not a variable. Each choice of $a$ gives rise to a different parametrization
of the evolution of the same time-dependent measure. System \eqref{achtung} shows
that to propagate such a measure all we require is the propagating velocity field
on the curve itself. 

\vspace{0.5cm}

Now let us return to the vortex sheets themselves. We assume vorticity
is of the form $\omega = \gamma \delta_{\ca{C}_t}$, and we parametrize
$\ca{C}_t$ by a function $y = y(\eta,t)$, with 
$\sigma = \gamma ds/d\eta$, as before.
The velocity associated to points $x$ {\it outside} the vortex sheet can be 
expressed  by the Biot-Savart law \eqref{eq:bsvart}: 
\begin{equation} \label{BS} 
u(x,t) = \int_{\eta_0}^{\eta_1} K(x-y(\eta,t))\sigma(\eta,t) d\eta, 
\end{equation}
with $K$ given by \eqref{K}, since $\sigma$ already includes the element of
length of the curve. 

It is a well-known fact that the flow \eqref{BS}, induced by the vortex sheet, is discontinuous 
across the sheet. More precisely, the normal component of $u$ at $\ca{C}_t$ is continuous, 
whereas the tangential component has a jump discontinuity with magnitude given precisely 
by $\gamma$. These are non-trivial facts and the reader may consult 
\cite{saffman} for a thorough discussion. 

Observe that the motion of the curve $\ca{C}_t$ is completely determined 
by the extension to $\ca{C}_t$ of the normal component of $u$. 
On the other hand, to propagate the density $\gamma$ in an 
explicit manner one needs to make a choice of tangential component of velocity on 
the sheet. This choice must take into account the nonlinear nature of the problem,
in a way that we will explore later.  For the present discussion, let us simply
consider the standard choice, which is to prescribe the velocity of the sheet 
as the {\it arithmetic mean} of the limit velocity from each side of the sheet. 
This arithmetic mean, when calculated using the velocity defined by \eqref{BS}, 
can be expressed as a principal value integral in the following way:

\begin{equation} \label{plemmelj}
\ca{U}[y;\sigma](\eta,t) \equiv p.v. \int_{a_0}^{a_1} K(y(\eta,t)-y(\eta^{\prime},
t))\sigma(\eta^{\prime},t) 
d\eta^{\prime} 
\end{equation}
\[ = \lim_{\varepsilon \to 0^+} \int_{|y(\eta,t) - y(\eta^{\prime},t)|\geq \varepsilon} 
K(y(\eta,t)-y(\eta^{\prime},t)) \sigma(\eta^{\prime},t) d\eta^{\prime}, \]
see \cite{saffman}
 
We use the vector field $\ca{U}$ to propagate $\omega = \gamma \delta_{\ca{C}_t}$.
We assume that the evolution of $\omega$ can be described by a system of the form 
\eqref{achtung} with transporting velocity $v=\ca{U}$ as follows

\begin{equation} \label{genbr}
\left\{ \begin{array}{l} y_t + a(\eta,t) y_{\eta} = \ca{U}[y;\sigma] \\ 
\sigma_t + (a(\eta,t)\sigma)_{\eta} = 0,
\end{array} \right. \end{equation}
 
System \eqref{genbr} is a general form of the classical Birkhoff-Rott system.
One may close system \eqref{genbr} by prescribing $a$. 
For instance, if one assumes that the vortex sheet is the graph of a 
function of $x$, and parametrizes it using $x$ itself, 
then $a$ is the first component of $\ca{U}[y;\sigma]$, see \cite{MP}. 
The choice of a {\it Lagrangian} parametrization, i.e., such
that $\partial y / \partial t = \ca{U}[y;\sigma]$ corresponds to 
choosing $a=0$. The scalar $a$ measures how much the evolution of a 
chosen parametrization fails to be Lagrangian. The circulation parametrization
$\Gamma = \Gamma(s,t) \equiv \int_0^s \gamma(s^{\prime},t)\,ds^{\prime}$, 
with $s$ being arclength, is a special case of Lagrangian parametrization for which 
$\sigma \equiv 1$. It gives rise to the traditional form of the Birkhoff-Rott
equations, 
\[ \partial_t \overline{z} = \frac{1}{2\pi} p. v. \int 
\frac{1}{z(\Gamma,t) - z(\Gamma^{\prime},t)} \, d\Gamma^{\prime},\]
where we have switched to complex variable notation for the position of the sheet. 

In this work, we will choose to parametrize vortex sheets by arclength. 
In this case, the function $a$ must become another unknown and the Birkhoff-Rott 
equations become:
\begin{equation} \label{bral}
\left\{ \begin{array}{l} \xi_t + a(s,t) \xi_s = \ca{U}[\xi;\gamma] \\ \gamma_t + (a(s,t)\gamma)_s = 0
\\ \left| \frac{\partial \xi}{\partial s} \right|= 1.
\end{array} \right. \end{equation}
One may also fix the origin of the arclength parametrization by taking $a(0,t)\equiv 0$, see \cite{wu}.

\vspace{0.5cm}

{\bf Remark:} How does one justify the use of the arithmetic mean in extending the Biot-Savart
velocity to the vortex sheet? This is, in a sense, the key issue behind the present work.
One could offer a convincing argument, approximating the evolution of the vortex sheet by
desingularizing the Biot-Savart kernel, using Lemma \ref{lemmachtung} for this situation to get
an approximate Birkhoff-Rott system and passing to the limit. This can be
done rigorously if we assume that the approximate evolution is convergent and 
if we also show that natural desingularizations indeed lead to the principal value integral 
\eqref{plemmelj}. A stronger version of the argument outlined above was carried out by 
Benedetto and Pulvirenti, who proved that the evolution of vortex sheets by Birkhoff-Rott is the
asymptotic description of thin shear bands under the Euler equations, see \cite{BP}.  
Another possibility is to argue that the arithmetic mean is the only extension that 
leads to vortex sheet evolution compatible with 2D Euler. This approach has been 
carried out in several manners, see \cite{birkhoff,majdabertozzi,MP,saffman}. 
In all cases, smoothness of the vortex sheet and its density have been assumed. One of the
motivations of the present work is to determine how irregular the vortex sheet can be,
while retaining the compatibility of the choice of arithmetic mean in \eqref{achtung} with 
incompressible 2D Euler.

\subsection{Weak form of the Birkhoff-Rott system}

Our goal is to compare solutions of the Birkhoff-Rott equations and of
the vorticity equation having limited smoothness.  To do so, we require a weak
formulation of the Birkhoff-Rott equations. Such a weak
formulation can be obtained by formally
substituting the vorticity \eqref{eq:sheetvort} into the traditional
weak formulation of the vorticity equation, given by \eqref{eq:twvort} 
and replacing $u$ by $\ca{U}$. We use the traditional rather than the modern 
weak formulation \eqref{wvfid} of the vorticity equation because, as we have 
discussed in the previous subsection, the hallmark of the Birkhoff-Rott equation 
is the choice of the arithmetic mean in extending velocity to the sheet, or 
equivalently, the introduction of the principal value in
the integral defining the velocity. The issue is that the 
velocity does not appear in the modern weak formulation of 
the vorticity equation.  
Plugging \eqref{eq:sheetvort} into
\eqref{eq:twvort} and replacing $u$ by $\ca{U}$ yields 
\begin{multline}
  \label{eq:wbr}
  \int_0^T\int_{s_0}^{s_1} \gamma(s,t) 
\lB \varphi_t(\xi(s,t),t)
+\ca{U}[\xi,\gamma](s,t) \cdot\grad\varphi(\xi(s,t),t)\rB\,ds\,dt
\\+\int_{s_0}^{s_1} \gamma_0(s)\varphi(\xi_0(s),0)\,ds=0.
\end{multline}

\begin{definition} \label{wbr} 
Let $\ca{C}_t = \{\xi = \xi(s,t)\mid s_0(t) \leq s \leq s_1(t)\}$ be a 
  rectifiable curve for each $t\in[0,T)$. Let $\gamma =
  \gamma(s,t) \in L^{\infty}((0,T);L^1(ds))$ be such that
  $\ca{U}[\xi;\gamma]$ is defined and $\gamma \,\ca{U}[\xi;\gamma] \in
  L^{\infty}((0,T);L^1(ds))$. We say that $(\gamma,\xi)$ is a weak
  solution of the Birkhoff-Rott equations with initial data
  $(\gamma_0,\xi_0)$ if \eqref{eq:wbr} holds for every test function $\varphi \in
  C^{\infty}_c([0,T) \times \real^2)$.
\end{definition}

 Our next step is to prove that Definition \ref{wbr} is compatible with
 the Birkhoff-Rott system \eqref{bral}. For vortex sheets of finite length 
 we will need to  supplement \eqref{bral} by (moving) boundary conditions.  
 It is natural to assume
that the endpoints $s_0(t)$ and $s_1(t)$ of the curve are material
trajectories, although this hypothesis will only be needed when the
vorticity density is nonzero there. We will therefore assume that
\begin{equation} \label{grumble}
\gamma(s_i(t),t)\lb s_i^{\prime}(t)-a(s_i(t),t) = 0\rb, \; i=0,1.
\end{equation}
Even when the vortex sheet has infinite length, if it is contained in
a bounded region then some form of boundary condition is still
needed. It will suffice to assume in this case that 
\begin{equation}
  \label{eq:gamlone}
  \int_{s_0(t)}^{s_1(t)} |\gamma(s,t)|\,ds\le c<\infty,
\end{equation}
which is reasonable since it makes the mass of vorticity locally
finite.

We now show that for smooth enough functions $\xi$ and $\gamma$, being 
a weak solution of the Birkhoff-Rott equations using arclength
parametrization is equivalent
to satisfying the classical Birkhoff-Rott system \eqref{bral}:
\begin{theorem} \label{consist}
  Let $\xi = \xi(s,t)$, $\gamma = \gamma(s,t)$ and $a = a(s,t)$ be
  solutions of \eqref{bral} in $C^1_b(\Omega)$, where
  $\Omega\eqdef\{(s,t)\,|\,s_0(t) \leq s \leq s_1(t), 0 \leq t <T\}$.
  Assume also that if an endpoint $s_j(t)$ is finite then
  \eqref{grumble} is satisfied, while if $s_j(t)$ is infinite but
  $\xi(s_j(t),t)$ is finite then \eqref{eq:gamlone} holds.  Then
  $(\gamma,\xi)$ is a weak solution of the Birkhoff-Rott equations
  with initial data $(\gamma(s,0),\xi(s,0))$ in the sense of
  Definition \ref{wbr}.
  
  Conversely, suppose that $(\gamma,\xi)$ is a weak solution of the
  Birkhoff-Rott equations with initial data $(\gamma_0,\xi_0)$ in the
  sense of Definition \ref{wbr}, where $s$ is an arclength
  parameter and at each fixed time $\xi$ is one-to-one except that
  $\xi(s_1(t),t)$ is allowed to equal $\xi(s_0(t),t)$.  Assume in
  addition the following prescribed regularity:
\begin{enumerate}
\item The parametrization $\xi$ and the density $\gamma$ are $C^1_b$ on $\Omega$;
\item The velocity $\ca{U}[\xi;\gamma]$ is  $C^0_b$ on $\Omega$;
\item If $|s_i| < \infty$ then $s_i$ is $C^1_b([0,T))$;
\end{enumerate}
Then $\gamma$ and $\xi$ satisfy \eqref{bral} and \eqref{grumble} with 
$a\eqdef \xi_s\cdot\lb \ca{U}-\xi_t\rb$, except that the equation for
$\xi$ need not hold in any open set on which $\gamma$ vanishes identically.
\end{theorem}

\vspace{0.5cm}

{\bf Remarks:}  
\begin{enumerate}
\item Suitably interpreted, Theorem 6.1.1 of \cite{MP} shows a version
  of the first half of Theorem~\ref{consist} for vortex sheets
  parametrized by one of the components of $x$.  To see this, note
  that, as remarked on the next page there, the integral $\int \omega
  u$ is not well-defined on the curve but must be given meaning via a
  principal-value integral.
\item The proviso that the evolution equation for $\xi$ need not hold
   where $\gamma$ vanishes is reasonable, since such regions are in
   essence not really part of the vortex sheet. Furthermore, the
   evolution equation for $\gamma$ implies that if it is nonzero
   everywhere on the sheet initially it will remain so at later times.
\end{enumerate}

\vspace{0.5cm}

\begin{proof}
  Let $\xi = \xi(s,t)$, $\gamma = \gamma(s,t)$ and $a = a(s,t)$ be
  $C^1_b$ solutions of \eqref{bral}.  In the following calculations we
  will assume that $(s_0(t),s_1(t))$ is a bounded interval for each
  $0\leq t < T$, but the case when either or both $s_j$ are infinite
  will also be treated.  Let $\varphi$ belong to $C^{\infty}_c([0,T)
  \times \real^2)$. Multiplying the equation for $\gamma$ in
  \eqref{genbr} by $-\varphi(\xi(s,t),t)$, integrating over $s$ and
  $t$, and then integrating by parts yields
  \begin{equation}\label{eq:intphigam}
    \begin{aligned}
      0&=-\int_0^T 
\int_{s_0(t)}^{s_1(t)} \lb \gamma_t+\lpp a\gamma\rpp_s\rb 
\varphi(\xi(s,t),t)\,ds\,dt
\\&= \int_0^T 
\int_{s_0(t)}^{s_1(t)} \gamma\lb \varphi_t+\xi_t\cdot\grad\varphi+
a\xi_s\cdot\grad\varphi\rb \,ds\,dt
\\&\fakeq+\int_{s_0(0)}^{s_1(0)} \gamma(s,0)\varphi(\xi(s,0),0)\,ds 
+\int_0^T \gamma\varphi\eval{s=s_j(t)} 
\lb s_j'(t)-a(s_j(t),t)\rb\Eval{j=0}{j=1}\,dt.
    \end{aligned}
  \end{equation}
If $s_j(t)$ is finite then the corresponding boundary term vanishes by
\eqref{grumble}. If $s_j(t)$ is infinite and $\xi(s_j(t),t)$ is also
infinite then that boundary term vanishes because $\varphi$ has
compact support. If $s_j(t)$ is infinite but $\xi(s_j(t))$ is finite
then the integrability of $\gamma$ combined with the boundedness of
its derivative implies that $\gamma$ tends to zero as $s\to s_j$, so
that the boundary term still vanishes. Upon taking into account the
equation for $\xi$ in \eqref{genbr}, \eqref{eq:intphigam} reduces to 
\begin{equation*}
\begin{aligned}
  0&=\int_0^T \int_{s_0(t)}^{s_1(t)} \gamma(s,t)
\lb \varphi_t(\xi(s,t),t) 
+ \ca{U}[\xi,\gamma](s,t)\cdot\grad\varphi(\xi(s,t),t)\rb\,ds\,dt
\\&\fakeq+\int_{s_0(0)}^{s_1(0)} \gamma(s,0)\varphi(\xi(s,0),0)\,ds,
\end{aligned}
\end{equation*}
which shows that \eqref{eq:wbr} holds with $\gamma_0(s)= \gamma(s,0)$
and $\xi_0(s)=\xi(s,0)$. Since $\gamma$ lies in $L^1$ by assumption,
and the conditions on $\xi$ and $a$ imply that $\ca{U}$ is bounded,
the other conditions of Definition \ref{wbr} are also satisfied.

Conversely, let $\gamma$ and $\xi$ be a weak solution of Birkhoff-Rott
in the sense of Definition \ref{wbr} satisfying the regularity
assumptions in the statement. Again, we will assume in our computation
that $s_0(t)$ and $s_1(t)$ are finite; if not then the boundary terms
vanish. Let $\varphi \in C^{\infty}_c([0,T) \times \real^2)$. Taking
$\varphi(x,t)=e^{-t/\varepsilon}\psi(x)$ in \eqref{eq:wbr}, letting
$\varepsilon\to0$, and using the assumed regularity shows that
$\xi_0(s)=\xi(s,0)$ and $\gamma_0(s)=\gamma(s,0)$.  Integration by
parts in \eqref{eq:wbr} therefore yields
\begin{equation}\label{eq:wtos1}
  \begin{aligned}
    0&=\int_0^T\int_{s_0(t)}^{s_1(t)} 
\gamma(s,t)\lb \varphi_t(\xi(s,t),t)+\ca{U}\cdot\grad\varphi(\xi(s,t),t)\rb
\,ds\,dt
\\&\fakeq+\int_{s_0(0)}^{s_1(0)}\gamma(s,0)\varphi(\xi(s,0),0)\,ds 
\\&=\int_0^T\int_{s_0(t)}^{s_1(t)} 
\gamma(s,t)\lb \ddt\varphi(\xi(s,t),t)+\lB\ca{U}-\xi_t\rB
\cdot\grad\varphi(\xi(s,t),t)\rb
\,ds\,dt
\\&\fakeq+\int_{s_0(0)}^{s_1(0)}\gamma(s,0)\varphi(\xi(s,0),0) \,ds
\\&=\int_0^T\int_{s_0(t)}^{s_1(t)} \gamma \lB\ca{U}-\xi_t\rB
\cdot\grad\varphi(\xi(s,t),t)-\gamma_t \varphi(\xi(s,t),t)\,ds\,dt
\\&\fakeq-\int_0^T\gamma\varphi\eval{s=s_j(t)}s_j'(t)\Eval{j=0}{j=1}\,dt
  \end{aligned}
\end{equation}
Since $\xi(s,t)$ is an arclength parametrization at each fixed time,
$|\xi_s|=1$, and hence $\ca{U}-\xi_t=a(s,t)\xi_s+b(s,t)\xi_s^\perp$
with $a=\xi_s\cdot \lb \ca{U}-\xi_t\rb$ being continuous. Substituting
this into \eqref{eq:wtos1} yields
\begin{equation}\label{eq:wtos2}
  \begin{aligned}
    0&=\int_0^T\int_{s_0(t)}^{s_1(t)} 
\lb a\gamma \partial_s \varphi(\xi(s,t),t)-\gamma_t
\varphi(\xi(s,t),t)\rb\,ds\,dt
\\&\fakeq+\int_0^T\int_{s_0(t)}^{s_1(t)} 
b\gamma \xi_s^\perp\cdot\grad\varphi(\xi(s,t),t)\,ds\,dt
-\int_0^T\gamma\varphi\eval{s=s_j(t)}s_j'(t)\Eval{j=0}{j=1}\,dt.
  \end{aligned}
\end{equation}
Since no time derivatives are applied to $\varphi$ in
\eqref{eq:wtos2}, taking $\varphi=\eta_\varepsilon(t)\psi(x)$ with
$\eta_\varepsilon(t)\to\delta(t-t_0)$ shows that for every $t\in(0,T)$,
\begin{equation}\label{eq:wtos3}
  \begin{aligned}
    0&=\int_{s_0(t)}^{s_1(t)} 
\lb a\gamma \partial_s \psi(\xi(s,t))-\gamma_t
\psi(\xi(s,t))\rb\,ds
\\&\fakeq+\int_{s_0(t)}^{s_1(t)} 
b\gamma \xi_s^\perp\cdot\grad\psi(\xi(s,t))\,ds
-\gamma\psi\eval{s=s_j(t)}s_j'(t)\Eval{j=0}{j=1}.
  \end{aligned}
\end{equation}
Now pick any $s_*\in (s_0(t),s_1(t))$, and take
\begin{equation*}
\psi(x)=\xi_s^\perp(s_*,t)\cdot\lb x-\xi(s_*,t)\rb
\,\frac{\eta\lpp\frac{|x-\xi(s_*,t)|^2}{\varepsilon^2}\rpp}{\varepsilon},
\end{equation*}
where $\eta$ is an even $C^\infty_c$ function. 
Since 
\begin{equation*}
\xi(s,t)-\xi(s_*,t)=\xi_s(s_*,t)\cdot(s-s_*)+o(s-s_*), 
\end{equation*}
we obtain the estimates
$|\xi(s,t)-\xi(s_*,t)|=O(|s-s_*|)$, 
\begin{equation*}
\begin{aligned}
  \psi(\xi(s,t))&=o(s-s_*)\,
\frac{\eta\lpp\frac{|x-\xi(s_*,t)|^2}{\varepsilon^2}\rpp}{\varepsilon},
\\
  \xi_s^\perp(s_*,t)\cdot\grad\psi(\xi(s,t))&=
\frac{\eta\lpp\frac{|x-\xi(s_*,t)|^2}{\varepsilon^2}\rpp}{\varepsilon}
+\frac{o((s-s_*)^2)}{\varepsilon^2}
\frac{\eta'\lpp\frac{|x-\xi(s_*,t)|^2}{\varepsilon^2}\rpp}{\varepsilon},
\\
\noalign{and}
 \xi_s(s_*,t)\cdot\grad\psi(\xi(s,t))&=
\frac{o((s-s_*)^2)}{\varepsilon^2}
\frac{\eta'\lpp\frac{|x-\xi(s_*,t)|^2}{\varepsilon^2}\rpp}{\varepsilon}.
\end{aligned}
\end{equation*}
Plugging these estimates into \eqref{eq:wtos3} and noting that the
boundary term there vanishes for $\varepsilon$ sufficiently small
yields 
$0=c\gamma(s_*,t)b(s_*,t)+o(1)$ for some nonzero $c$, which shows that
$b$ times $\gamma$ vanishes identically.  Hence  \eqref{eq:wtos3}
reduces to 
\begin{equation}\label{eq:wtos4}
  \begin{aligned}
    0&=\int_{s_0(t)}^{s_1(t)} 
\lb a\gamma \partial_s \psi(\xi(s,t))-\gamma_t
\psi(\xi(s,t))\rb\,ds
-\gamma\psi\eval{s=s_j(t)}s_j'(t)\Eval{j=0}{j=1}.
  \end{aligned}
\end{equation}
In particular, 
\begin{equation}\label{eq:wtos5}
  \begin{aligned}
    0&=\int_{s_0(t)}^{s_1(t)} 
\lb a\gamma \partial_s \psi(\xi(s,t))-\gamma_t
\psi(\xi(s,t))\rb\,ds
  \end{aligned}
\end{equation}
for every $\psi$ that vanishes at the endpoints $s_j(t)$.  This
implies that $a\gamma$ is differentiable with respect to $s$ for
$s_0(t)<s<s_1(t)$ and
$(a\gamma)_s=-\gamma_t$, i.e., the equation for $\gamma$ in
\eqref{genbr} holds. Furthermore, since by construction
$\xi_t+a\xi_s-\ca{U}=b\xi_s^\perp$, the fact that $b\gamma$ vanishes
identically shows that the evolution equation for $\xi$ holds
wherever $\gamma$ is nonzero. Since the expression
$\xi_t+a\xi_s-\ca{U}$ is continuous it must then vanish except on open
sets where $\gamma$ vanishes identically. 

Finally, when at least one $s_j$ is finite then
integrating by parts in \eqref{eq:wtos4} now shows that
\begin{equation*}
  0=\gamma\psi\eval{s=s_j(t)}\lb s_j'(t)-a(s_j(t),t)\rb\Eval{j=0}{j=1},
\end{equation*}
and this implies \eqref{grumble}
whether the $s_j$ are distinct or not.
\end{proof}

\vspace{0.5cm}

{\bf Remark:}  One of the issues that made 
the proof above long was our concern in including as many plausible
examples of vortex sheet evolution as possible. For example, the Kaden
spiral, periodic sheets and closed sheets all may be considered as
smooth solutions. 
 
\section{Equivalence of weak formulations}

This section is also divided into three subsections. In the first one we
recall the concept of regular curves and G. David's result on Cauchy integrals
on regular curves. In the second subsection we state and prove our main 
result, the criterion for equivalence of Birkhoff-Rott and 2D Euler for 
vortex sheet evolution. In the final subsection we recall an example that
establishes the sharpness of the criterion presented.

\subsection{Regular Curves}

Let us begin by recalling the concept of {\it regular curve}. A
rectifiable curve $\ca{C}$ is called {\it regular} if there exists a
constant $A>0$ such that for any disk $D_r$ of radius $r>0$,
\begin{equation} \label{curveconst}
|\ca{C} \cap D_r| \leq Ar,
\end{equation} 
where $|\ca{C} \cap D_r|$ represents the length of this segment of curve. 

Let $\ca{C} \equiv \{ \xi = \xi(s) \}$ be a rectifiable curve
parametrized by arclength, with $s \in (a_0,a_1)$. Let $\ca{U}_{\ast}$
be the maximal operator associated to $\ca{U}$, i.e.,
\[ \ca{U}_{\ast}[\gamma] \equiv \sup_{\varepsilon > 0} \left| \int_{|\xi(s) - 
\xi(s^{\prime})|\geq \varepsilon} 
K(\xi(s)-\xi(s^{\prime})) \gamma(s^{\prime}) ds^{\prime} \right|. \]

\begin{theorem} \label{david} (G. David, \cite{david})
  Suppose that $\ca{C}$ is a regular curve and let $1<p<\infty$. Then
  the maximal operator $\gamma \mapsto \ca{U}_{\ast}[\gamma]$ is a
  bounded sublinear operator from $L^p(ds)$ into $L^p(ds)$.
  Conversely, if there exists a continuous linear operator $\ca{U}:
  L^2(ds) \to L^2(ds)$ such that, for any $\gamma \in C^0_c(ds)$ and
  for any $s_0$ such that $\xi(s_0)$ does not belong to
  $\xi(\supp(\gamma))$, it holds that $\ca{U}[\xi;\gamma](s_0) = \int
  K(\xi(s_0)-\xi(s^{\prime})) \gamma(s^{\prime}) ds^{\prime}$, then
  $\ca{C}$ is a regular curve.
\end{theorem}
 
\vspace{0.5cm}

{\bf Remark:}
\begin{enumerate}
\item The result above was originally stated as a
characterization of the rectifiable curves $\ca{C}$ in the complex
plane such that the Cauchy integral defines a bounded operator from
$L^2(\ca{C})$ to itself, namely, Theorem 2 in \cite{david}. One may
identify the integral in the definition of $\ca{U}_{\ast}$ with a
Cauchy integral by introducing the usual identification of $\real^2$
with $\complex$. In addition, we stated David's result in $L^p$,
thereby incorporating the comment made immediately after the proof on
page 174 of \cite{david}.

\item As remarked in \cite{david}, it follows that in the
case of parametrization by arclength, $\ca{U}[\xi;\gamma]$ is defined
$ds$-almost everywhere. Furthermore, $\ca{U}$ defines a continuous
linear operator from $L^p(ds)$ to itself.
\end{enumerate}

\subsection{The Equivalence Theorem}

We are now ready to state and prove our main result.

\begin{theorem} \label{kabum}
  Let $\ca{C}_t = \{\xi = \xi(s,t)\,|\, s_0(t) \leq s \leq a_1(t)\}$
  be a regular curve parametrized by arclength, $0\leq t<T$. Assume
  that the constant $A$ in the definition of regular curve,
  \eqref{curveconst}, may be chosen independently of $t$. Let $\gamma
  \in L^{\infty}([0,T);L^2(ds)\cap L^1(ds))$. Then $(\gamma,\xi)$ is a
  weak solution of the Birkhoff-Rott equations with initial data
  $(\gamma_0,\xi_0)$ if and only if $\omega = \gamma
  \delta_{\ca{C}_t}$ is a weak solution of the weak vorticity
  formulation with initial data $\omega_0= \gamma_0\delta_{\ca{C}_0}$,
  $\ca{C}_0 = \{\xi=\xi_0(s)\}$.
\end{theorem}

Before we proceed with the proof, let us emphasize that the {\it only if}
portion of this result assumes that the weak solution of 2D Euler has the
structure $\gamma \delta_{\ca{C}_t}$. Although existence of weak solutions
with initial data of this kind has been established in certain cases, 
their structure is not known {\it a priori}.

\begin{proof}
  There are two steps in this proof. The first step is to show that
  the identities \eqref{eq:wbr} and \eqref{vsws} are the same under the
  hypothesis' of this theorem. The second is to show that the
  regularity requirements in Definitions \ref{wbr} and \ref{wvf} are
  equivalent in this case as well.
  
  We begin by showing that the identities \eqref{eq:wbr} and
  \eqref{vsws}, involving the test function $\varphi \in
  C^{\infty}_c([0,T)\times \real^2)$, are the same. First note that we
  need only consider the nonlinear term in each identity. We will show
  that, under our hypothesis', we have
\[ \int_0^T \int_{s_0(t)}^{s_1(t)}   
\nabla \varphi(\xi(s,t),t) \cdot \ca{U}[\xi;\gamma](s,t) \, \gamma(s,t) dsdt  \]
\[= \int_0^{T} \int_{s_0(t)}^{s_1(t)} \int_{s_0(t)}^{s_1(t)} 
H_{\varphi}(\xi(r,t),\xi(s,t),t)\gamma(r,t)\gamma(s,t) \, drdsdt, \]
for any test function $\varphi$. We start by recalling that 
\[\ca{U}[\xi;\gamma] = \ca{U}[\xi;\gamma](s,t) 
= \lim_{\varepsilon \to 0^+} \int_{|\xi(s,t) - \xi(r,t)|
\geq \varepsilon} K(\xi(s,t)-\xi(r,t)) \gamma(r,t) \,dr. \]

Fix $\varepsilon>0$. Denote by $\Delta_{\xi}(r,s,t) \eqdef \xi(s,t) - \xi(r,t)$. Note that
\[ \int_0^T \int_{s_0(t)}^{s_1(t)}   \nabla \varphi(\xi(s,t),t) \cdot  
\left(\int_{|\Delta_{\xi}(r,s,t)|\geq 
\varepsilon} K(\Delta_{\xi}(r,s,t)) \gamma(r,t) \,dr\right)
\gamma(s,t) dsdt  \]
\begin{equation} \label{clearly} 
= \int_0^T \int_{s_0(t)}^{s_1(t)}   \int_{|\Delta_{\xi}(r,s,t)|\geq \varepsilon} 
\nabla \varphi(\xi(s,t),t) \cdot   K(\Delta_{\xi}(r,s,t)) 
\gamma(r,t)   \gamma(s,t) \,dr dsdt
\end{equation}
\[=- \int_0^T \int_{s_0(t)}^{s_1(t)} 
\int_{|\Delta_{\xi}(r,s,t)|\geq \varepsilon}  \nabla \varphi(\xi(r,t),t) \cdot  
 K(\Delta_{\xi}(r,s,t)) \gamma(s,t) 
\gamma(r,t)\,ds  drdt,  \]
exchanging $s$ with $r$ and using the antisymmetry of the kernel $K$,  
\begin{equation} \label{blic1}
 = - \int_0^T \int_{s_0(t)}^{s_1(t)}  
\int_{|\Delta_{\xi}(r,s,t)|\geq \varepsilon} 
\nabla \varphi(\xi(r,t),t) \cdot  K(\Delta_{\xi}(r,s,t)) \gamma(s,t) 
\gamma(r,t) dr \,dsdt,  
\end{equation}
using Fubini's theorem.  

Hence, adding $1/2$ of \eqref{clearly} to $1/2$ of \eqref{blic1}, we find that 
\[ \int_0^T \int_{s_0(t)}^{s_1(t)}   \nabla \varphi(\xi(s,t),t) \cdot  
\left(\int_{|\Delta_{\xi}(r,s,t)|\geq \varepsilon} 
K(\Delta_{\xi}(r,s,t)) \gamma(r,t) \,dr\right) \gamma(s,t) dsdt  \]
\[= \int_0^{T} \int_{s_0(t)}^{s_1(t)} 
\int_{|\Delta_{\xi}(r,s,t)|\geq \varepsilon}  
H_{\varphi}(\xi(r,t),\xi(s,t),t)\gamma(r,t)\gamma(s,t) \, drdsdt. \]

The right hand side of the identity above converges to the nonlinear
term in \eqref{vsws} as $\varepsilon \to 0$, by the Dominated
Convergence Theorem, since $H_{\varphi}$ is bounded and $\gamma$ was
assumed to be integrable.

To establish that \eqref{eq:wbr} and \eqref{vsws} are the same it
remains to show that
\[ \int_0^T \int_{s_0(t)}^{s_1(t)}   
\nabla \varphi(\xi(s,t),t) \cdot \ca{U}[\xi;\gamma](s,t) \, \gamma(s,t) dsdt  \]
\[ = \lim_{\varepsilon \to 0} \int_0^T \int_{s_0(t)}^{s_1(t)}   
\nabla \varphi(\xi(s,t),t) \cdot  
\left(\int_{|\Delta_{\xi}(r,s,t)|\geq \varepsilon} 
K(\Delta_{\xi}(r,s,t)) \gamma(r,t) \,dr\right)
 \, \gamma(s,t) dsdt.  \] 

We first note that 
\[\left|\int_{|\Delta_{\xi}(r,s,t)|\geq \varepsilon} 
K(\Delta_{\xi}(r,s,t)) \gamma(r,t) \,dr \right| \leq 
\ca{U}_{\ast}[\gamma],\]
where $\ca{U}_{\ast}$ is the maximal operator associated 
to $\ca{U}[\cdot,\xi]$ introduced in Section 2.  
Since $\gamma \in L^{\infty}((0,T);L^2(ds))$, 
it follows from Theorem \ref{david} that, $\ca{U}_{\ast}[\gamma]
\in L^2(ds)$, for almost all time. Furthermore, since $\ca{U}_{\ast}$ is a bounded
sublinear operator from $L^2$ to itself, we have in fact that
$\ca{U}_{\ast}[\gamma] \in L^{\infty}((0,T);L^2(ds))$. Therefore
$\nabla \varphi \cdot \ca{U}_{\ast}[\gamma] \gamma \in
L^{\infty}((0,T);L^1(ds))$. The desired identity follows by dominated
convergence.

The second step in this proof is to examine regularity conditions in
Definition \ref{wbr} and Definition \ref{wvf}.  First note that any
density $\gamma$ and curve $\xi$ satisfying the hypothesis of this
Theorem will satisfy the regularity requirements of Definition
\ref{wbr}. Indeed, by Theorem \ref{david}, $\ca{U}[\xi;\gamma] \in
L^{\infty}((0,T);L^2(ds))$ so that $\gamma \, \ca{U}[\xi;\gamma]\in
L^{\infty}((0,T);L^1(ds))$. On the other hand, the fact that
$\gamma \in L^{\infty}((0,T);L^1(ds))$ implies that $\omega\in
L^\infty((0,T);\ca{BM}_c)$.

Under the conditions of the theorem, the vorticity also
  automatically satisfies the condition 
\begin{equation}\label{eq:c2}
\omega \in
  L^{\infty}((0,T);H^{-1}_{\loc}(\real^2))
\end{equation}
needed to obtain the
  equivalence of the weak vorticity formulation with the definitive
  weak velocity formulation. In addition, $\omega$ also satisfies the
  condition 
\begin{equation}\label{eq:c3}
\omega \in Lip((0,T);H^{-L}_{\loc}(\real^2))\quad\text{for some
  $L > 1$}
\end{equation}
that ensures the compactness of bounded sequences of solutions.

To prove \eqref{eq:c2}, we begin by observing that, given any disk
$D_r$, of radius $r>0$ in the plane we have
\[ |\omega(D_r)| = \left| \int_{\ca{C}_t \cap D_r} \gamma ds \right| 
\leq \|\gamma\|_{L^2}^{1/2} (Ar)^{1/2} \leq Cr^{1/2},\] by
Cauchy-Schwarz and using the definition of regular curve, for some
constant $C>0$ independent of time.  Hence the vorticity
$\omega$ belongs to $L^{\infty}((0,T);\widetilde{M}^{4/3}(\real^2))$,
where $\widetilde{M}^{4/3}$ is a Morrey space of
measures~\cite{lopes2}.  It was shown in \cite{lopes2}, Theorem 4.3,
that for any $p > n/2$, $\widetilde{M}^p (\real^n) \cap
\ca{BM}_{\loc}(\real^n)$ is (compactly) contained in
$H^{-1}_{\loc}(\real^n)$.

Now consider \eqref{eq:c3} We have already shown that, since
$(\gamma,\xi)$ was assumed to satisfy Definition \ref{wbr}, then
$(\gamma,\xi)$ satisfies identity \eqref{vsws} for any test function
$\varphi \in C^{\infty}_c([0,T) \times \real^2)$. In particular, if
$\varphi \in C^{\infty}_c((0,T) \times \real^2)$ then
\[\int_0^T\int_{s_0(t)}^{s_1(t)} \varphi_t(\xi(s,t),t)\gamma(s,t) \, dsdt \]
\[= - \int_0^{T} \int_{s_0(t)}^{s_1(t)} 
\int_{s_0(t)}^{s_1(t)} H_{\varphi}(\xi(r,t),\xi(s,t),t)
\gamma(r,t)\gamma(s,t) \, drdsdt. \]
Observe that $\|H_{\varphi}(\cdot,\cdot,t)\|_{L^{\infty}(\real^4)} 
\leq \|\varphi(\cdot,t)\|_{W^{2,\infty}(\real^2)}$. Note also that
\[\int_0^T\int_{s_0(t)}^{s_1(t)} \varphi_t(\xi(s,t),t)\gamma(s,t) \,
dsdt = \int_0^T\int_{\real^2} \varphi_t \omega (x,t) \,dxdt.\] 
In view of these observations, together with the hypothesis' on
$\gamma$, we deduce that
\[ \left| \int_0^T\int_{\real^2} \varphi_t \omega \,dxdt \right|
 \leq C\| \varphi \|_{L^1((0,T);W^{2,\infty}(\real^2))},\]
for any test function in $C^{\infty}_c((0,T) \times \real^2)$,
\[\leq C\| \varphi \|_{L^{1}((0,T);H^4(\real^2))},\]
by the Morrey inequality. By duality this implies that the distribution 
$\omega_t \in L^{\infty}((0,T);H^{-4}(\real^2))$, which in turn 
gives \eqref{eq:c3} with $L=4$. 

\end{proof}

{\bf Remark:} The hypothesis
that $\gamma \in L^{\infty}((0,T);L^1(ds))$ is redundant, as it is
explicitly present in Definition \ref{wbr} and implicitly present in
Definition \ref{wvf}. We included this hypothesis in the statement of
the theorem merely for the sake of clarity.

\subsection{Sharpness}
Next we will observe that the condition $\gamma \in L^2(ds)$, needed
above to prove the equivalence of the weak formulations of 2D Euler
and Birkhoff-Rott, is sharp.  To see this, we recall the classical
example of the Prandtl-Munk vortex sheet, see \cite{munk19}, also known
as the elliptically loaded wing, see \cite{krasny87}.  We consider flow
generated by the initial vortex sheet given by
\[\ca{C}_0 = \{\xi_0(s) = (s,0), -1 \leq s \leq 1\} 
\,\mbox{ and } \,\gamma_0(s) = \frac{s}{\sqrt{1-s^2}}. \] By complex
variable methods it can be shown that the velocity
$\ca{U}_0[\xi_0;\gamma_0](s) \equiv (0,-1/2)$, see \cite{saffman},
section 6.2. This is pointwise true for $-1<s<1$, but unclear at $s =
\pm 1$, since the flow velocity is infinite near the tips of the
sheet.

The fact that the arithmetic mean $\mathcal{U}$ is constant 
suggests that $\xi(s,t) = (s, -t/2)$, $\gamma(s,t) = \gamma_0(s)$
describes the evolution of this vortex sheet.  In fact, $(\gamma,\xi)$
is indeed a weak solution of the Birkhoff-Rott equations with initial
data $(\gamma_0,\xi_0)$ in the sense of Definition \ref{wbr}. The
verification of this fact is straightforward. On the other hand, it
was shown in \cite{lopes1} that $\omega = \gamma \delta_{\ca{C}_t}$
is not a weak solution of the incompressible 2D Euler equations with
initial data $\gamma_0 \delta_{\ca{C}_0}$ in the sense of Definition
\ref{wvf}.  We observe that $\gamma \in L^1(ds) \cap L^p(ds)$ for all
$p<2$ but not for $p=2$, which shows the sharpness of the $L^2$
condition on $\gamma$ in Theorem \ref{kabum}.

\section{Final Remarks}

Recently, G. Lebeau proved that any solution of the Birkhoff-Rott system 
consisting of a closed vortex sheet which is $C^{1,\alpha}$ at a time $t_0$
is real-analytic for $t \neq t_0$, see \cite{lebeau}. This is a consequence of the 
elliptic nature of the Birkhoff-Rott equations. In  \cite{wu}, Sijue Wu announced 
an improvement of Lebeau's result, in which the conclusion is mantained if the vortex 
sheet is required to be a chord-arc curve and satisfies an additional technical 
condition. Thus, one might conclude that there are no irregular solutions 
of Birkhoff-Rott and wonder whether there is any reason to formulate a general 
theory of equivalence of Euler and weak Birkhoff-Rott. In that regard, we observe 
that the Prandtl-Munk vortex sheet is an example of an irregular solution of
Birkhoff-Rott which is not covered by Lebeau's work because it is not a closed
sheet, so that the context behind our result is not entirely empty. We also note 
that the class of regular curves is broader than chord-arc, hence solutions which are 
regular curves but which are not real analytic are plausible. In particular, the 
logarithmic spirals of Prandtl and Kambe fit into our discussion. 

In Theorem \ref{kabum} on the equivalence between the Euler and Birkhoff-Rott descriptions of 
vortex sheet motion we restrict our attention solely to regular curves. This hypothesis
is motivated by hindsight, since, David's result includes the fact that the 
Cauchy integral on a rectifiable curve is a bounded operator on $L^2$ only if the curve is 
regular, as stated in Theorem \ref{david}. Although the Birkhoff-Rott 
equations might make sense even if the principal-value integral does not give rise to
a bounded operator in $L^2$, this becomes a substantial complication. However, the 
hypothesis that the sheet be regular is a big limitation in our work, since 
the conjectured behavior past singularity formation for periodic sheets is 
the development of infinite length algebraic spirals, which are hence not regular.

Let us recall that the existence of vortex sheet evolution, at least with distinguished sign 
density, has been established from the implicit approach. Why then would we still 
want to solve Birkhoff-Rott? The main motivation is that the Birkhoff-Rott equations 
carry a much more precise description of the flow than the weak formulation of 2D Euler. 
In particular, the Birkhoff-Rott equations assume {\it a priori} that the evolution of
a vortex sheet retains a curve-like structure, whereas the structure of the weak solutions 
given, for instance by Delort's Theorem, is not known. One interesting open problem is to
prove that the support of a solution given by Delort's Theorem, with a smooth vortex sheet
as initial data, has Lebesgue measure zero (better yet, Hausdorff dimension less than 2) 
at almost every time.

Another interesting question is to find a meaningful way in which curves which are not regular 
may be regarded as solutions of Birkhoff-Rott. At this point it would be helpful to have 
examples of physically meaningful vortex sheets having infinite length algebraic spirals,
something which is not available at present. In this context it is worth mentioning a well-known 
open problem in the field, that of establishing existence of the Kaden spiral, or, more
generally, whether there exist self-similar, algebraic spiral solutions of the Birkhoff-Rott equations.        

\vspace{0.5cm}

\small{ {\it Acknowledgments:} The authors wish to thank Max Souza and Sijue Wu for helpful 
discussions. M. C. L. F.'s research was supported in part by CNPq grant \# 300.962/91-6. 
H. J. N. L.'s research was supported in part by CNPq grant \# 300.158/93-9.}

%

\end{document}